\def\Bbb#1{{\bf #1}}

\def\fnote#1{\footnote}
\def\blacksquare{\hbox{\vrule width 4pt height 4pt depth 0pt}}

\def\cwdash{\relbar\joinrel}
\def\cwleftpar#1#2{\leftskip #1 \rightskip #2 plus 1fill}
\def\cwrightpar#1#2{\leftskip #1 plus 1fill \rightskip #2}
\def\cwcenterpar#1#2{\leftskip #1 plus 1fill \rightskip #2 plus 1fill}
\def\cwfullpar#1#2{\leftskip#1\rightskip#2}

\def\cwoutdent#1#2{\llap{\hbox to #1{#2 \hss}}\ignorespaces}
\def\cwparbegin#1#2#3#4#5{
	\ifcase #1 \cwleftpar{#2}{#3}
	\or \cwrightpar{#2}{#3}
	\or \cwcenterpar{#2}{#3}
	\else \cwfullpar{#2}{#3}\fi
	\ifcase #4 \baselineskip = 1.5\baselineskip
	\or \baselineskip = 2\baselineskip
	\or \baselineskip = 3\baselineskip
	\else \baselineskip = 1\baselineskip\fi
	\ifdim #5 > 0in \else \noindent \fi
	\noindent\ignorespaces}
\documentclass[11pt]{article}
\begin{document}
%------------------------------------------------------------------------
%                               ChiWriter HEADER
%------------------------------------------------------------------------
\advance \vsize by -2\baselineskip
\def\makeheadline{
{\noindent \folio \par
\vskip \baselineskip }}
%------------------------------------------------------------------------
%                               Your Document
%------------------------------------------------------------------------

\vspace{2ex}
\noindent {\Huge Special bases for derivations of\\[0.3ex]
			tensor algebras}\\[0.4ex]
\noindent {\Large I. Cases in a neighborhood and at a point}

\vspace{2ex}

\noindent Bozhidar Z. Iliev
\fnote{0}{\noindent $^{\hbox{}}$Permanent address:
Laboratory of Mathematical Modeling in Physics,
Institute for Nuclear Research and \mbox{Nuclear} Energy,
Bulgarian Academy of Sciences,
Boul.\ Tzarigradsko chauss\'ee~72, 1784 Sofia, Bulgaria\\
\indent E-mail address: bozho@inrne.bas.bg\\
\indent URL: http://theo.inrne.bas.bg/$^\sim$bozho/}

\vspace{3ex}

{\bf \noindent Published: Communication JINR, E5-92-507, Dubna, 1992}\\[1ex]
\hphantom{\bf Published: }
http://www.arXiv.org e-Print archive No.~math.DG/0303373\\[2ex]

2003 PACS numbers: 02.40.Ma, 02.40.Vh,  04.20.Cv, 04.90.+h\\[2ex]

\noindent
{\small
The \LaTeXe\ source file of this paper was produced by converting a
ChiWriter 3.16 source file into
ChiWriter 4.0 file and then converting the latter file into a
\LaTeX\ 2.09 source file, which was manually edited for correcting numerous
errors and for improving the appearance of the text.  As a result of this
procedure, some errors in the text may exist.
}\\[0.2ex]

	\begin{abstract}
Necessary and sufficient conditions are investigated for the existence of
local bases in which the components of derivations of tensor algebras over
differentiable manifold vanish in a neighborhood or only at a single point.
The problem when these bases are holonomic or anholonomic is considered.
Attention is paid to the case of linear connections. Relations of these
problems with the equivalence principle are pointed out.
	\end{abstract}\vspace{2ex}

  {\bf I. INTRODUCTION}
\nopagebreak

\medskip
	In the theory of linear connections $[1, 2]$ problems connected with
existence of (local) bases in which the connection's components vanish at a
point $[2-7]$, along a curve $[3, 6]$ or in a neighborhood $[3, 5, 6]$ have
been considered. But with a very rare exceptions (see e.g. [7]) in the
literature only the torsion free case has been investigated. The aim of the
present work is generalization of these problems and their results to the
case of arbitrary derivations of the tensor algebra over a differentiable
manifold (see [2] or Sec. II), the curvature and torsion of which (as defined
below in Sec. II) are not a priori restricted somehow.

Mathematically the main purpose of this work is necessary and sufficient
conditions to be found for the existence of local (holonomic or anholonomic
[3]) bases (coordinates) in which the components of some derivation (of the
tensor algebra over a manifold) must vanish. If such bases exist, we
investigate the problem when they (or part of them) are holonomic.

 Physically the goal of the paper is to be shown that in gravity theories,
based, first of all, on linear connections, the equivalence principle is
identically satisfied because of their underlying mathematical structure.

The work is organized as follows. In Sec. II some notations and definitions
are introduced. Sec. III deals with the above pointed problems in a
neighborhood and  Sec. IV investigates them at a single point. In Sec. $V$
the connection of the considered mathematical problems with the equivalence
principle is shown. Sec. VI contains some concluding remarks.

\medskip
\medskip
  {\bf II. DERIVATIONS, THEIR COMPONENTS, CURVATURE AND TORSION}

\medskip
Let $D$ be a derivation of the tensor algebra over a given manifold $M [1,
2]$. By proposition 3.3 from $ch. I$ from [2] there exists unique vector
field $X$ and unique tensor field $S$ of type (1,1) such that $D=L_{X}+S$,
where $L_{X}$ is the Lie derivative along $X [1, 2]$ and $S$ is considered as
a derivation of the tensor algebra over $M [2]$.

If $S$ is a map from the set of $C^{1}$vector fields into the tensor fields
of type (1,1) and $S:X\mapstochar \cwdash \to S_{X}$, then the equation
\[
D^{S}_{X}=L_{X}+S_{X}\qquad (1)
\]
 defines a derivation of
the tensor algebra over $M$ for any $C^{1}$vector field $X [2]$. As the map
$S$ will hereafter be assumed fixed, such a derivation will be called an
$S$-derivation along $X$ and will be denoted by $D_{X}$. An $S$-derivation is
a map $D$ such that $D:X\mapstochar \to D_{X}$, where $D_{X}$is an
$S$-derivation along X.

Let $\{E_{i}, i=1,\ldots  n:=\dim(M)\}$ be a (coordinate or not $[3, 4])$
local basis of vector fields in the tangent to $M$ bundle. It is holonomic
(anholonomic) if the vectors $E_{1}, \ldots  , E_{n}$commute (don't commute)
$[3, 4]$. Let $T$ be a $C^{1}$ tensor field of type $(p,q), p$ and $q$ being
integers or zero(s), with local components $T^{i_{1}}_{j_{1}}$with respect to
the tensor basis associated to $\{E_{i}\}$. Here and below all latin indices,
may be with some subscripts, run from 1 to $n:= =\dim(M)$. Using the explicit
action of $L_{X}$ and $S_{X}$ on tensor fields [2] and the usual summation
rule on repeated on different levels indices, we find the components of
$D_{X}T$ to be
\[
(D_{X}T)^{i_{1}\ldots i_p}_{j_{1}\ldots j_q}
= X \bigl( T^{i_{1}\ldots i_p}_{j_{1}\ldots j_q} \bigr)
+\sum_{a=1}^{p} (W_{X})^{i_{a}}_{..k}
T^{i_{1}\ldots i_{a-1}ki_{a+1}\ldots i_p}_{j_{1}\ldots j_q}
\]
\[
+\sum_{b=1}^{q} (W_{X})^{k}_{.j_{b}}
T^{i_{1}\ldots i_p}_{j_{1}\ldots j_{b-1}kj_{b+1}\ldots j_q}
\qquad (2)
\]
where with $X(f)$ we denote the action of $X=X^{k}E_{k}$ on the $C^{1}$ scalar
function $f$, i.e. $X(f):=X^{k}E_{k}(f)$, and
\[
(W_{X})^{i}_{.j}:=(S_{X})^{i}_{.j}-E_{j}(X^{i})+C^{i}_{.kj}X^{k}\qquad (3)
\]
in which $C^{i}_{.kj}$defines the commutators of the basic vectors, i.e.

\[  [E_{j},E_{k}]=C^{i}_{.jk}E_{i}.\qquad (4)
\]
 By definition we shall call $(W_{X})^{i}_{.j}$ the components of $D_{X}$ as
they play with respect to $D_{X}$ the same role as the components of a linear
connection with respect to itself. In fact, from (2) we see that $(cf. (7)$
below)
\[
  D_{X}(E_{j})=(W_{X})^{i}_{.j}E_{i}.\qquad (5)
\]
 If we make a
change $\{E_{i}\}  \to \{E_{i^\prime }:=A^{i}_{i^\prime }E_{i}\}$, where
$A:= [ A^{i}_{i^\prime }]:= [ A^{i^\prime }_{i}]^{-1}$ is a nondegenerate
matrix function, then from (5) we can see that $(W_{X})^{i}_{.j}$ transform
into
\[
(W_{X})^{i^\prime }_{..j^\prime }=A^{i^\prime }_{i}A^{j}_{j^\prime
}(W_{X})^{i}_{.j}+A^{i^\prime }_{i}X(A^{i}_{j^\prime }),\qquad (6)
\]
which, if we introduce the matrices $W_{X}:=  (W_{X})^{i}_{.j}  $ and $W$
:$= := (W_{X})^{i^\prime }_{..j^\prime }  $, will read
\[ W=A^{-1}W_{X}A+X(A)^{  }_{  },\qquad (6^\prime )
\]
where the superscript is understand as a first matrix index  and
$X(A):=X^{k}E_{k}(A)= X^{k}E_{k}(A^{i}_{i^\prime })  $.

If $\nabla $ is a linear connection with local components $\Gamma
^{i}_{.jk}($see, e.g., $[1-3])$, then it is fulfilled
\[
\nabla _{X}(E_{j})=(\Gamma ^{i}_{.jk}X^{k})E_{i}.\qquad (7)
\]
Hence, comparing (5) and (7), we see that $D_{X}$is a covariant
differentiation along $X$ iff
\[
(W_{X})^{i}_{.j}=\Gamma ^{i}_{.jk}X^{k},\qquad (8)
\] for some functions
$\Gamma ^{i}_{.jk}$. Due to (3) or $(1), a$ linear connection $\nabla $ is
characterized by the map $S:X\mapstochar \to S_{X}$ such that
\[
S_{X}=\Sigma _{X}, \Sigma _{X}(Y):=\nabla _{X}(Y)-[X,Y],\qquad (9)
\]
where $[X,Y]=L_{X}Y$ is the commutator of the vector fields $X$ and $Y [2]$.

Let $D$ be an $S$-derivation and $X, Y$ and $Z$ be vector fields. The
curvature operator $R^{D}$and the torsion operator $T^{D}$ of $D$ are
\[
R^{D}(X,Y):=D_{X}D_{Y}-D_{Y}D_{X}-D_{[X,Y]},\qquad (10)
\]
\[
T^{D}(X,Y):=D_{X}Y-D_{Y}X-[X,Y],\qquad (11)
\]
which in a case of linear
connection reduce to the corresponding classical tensorial quantities (see
$[2, 3]$ and below (16) and (17)). The $S$-derivation $D$ will be called flat
(=curvature free) or torsion free if $R^{D}=0$ or $T^{D}=0$ respectively
$(cf. [2])$.

  If we use the representation (1) and $L_{X}Y=[X,Y]$, we get
\[
(R^{D}(X,Y))Z
=\{S_{X}S_{Y}-S_{Y}S_{X}+L_{X}(S_{Y})-L_{Y}(S_{X})\}Z+
\]
\[
\qquad
+ S_{Y}[X,Z]-S_{X}[Y,Z]-S_{[X,Y]}Z \qquad (12)
\]
\[
 T^{D}(X,Y)=S_{X}Y-S_{Y}X+[X,Y].\qquad (13)
\]
Analogously, by means of (2), we find the local expressions:
\[
(R^{D}(X,Y))^{i}_{.k}
=X((W_{Y})^{i}_{.k})-Y((W_{X})^{i}_{.k}) +(W_{X})^{i}_{.l}(W_{Y})^{l}_{.k} -
(W_{Y})^{i}_{.l}(W_{X})^{l}_{.k}-
\]
\[
\qquad -(W_{[X,Y]})^{i}_{.k},\qquad (14)
\]
\[
(T^{D}(X,Y))^{i}
=(W_{X})^{i}_{.l}Y^{l}-(W_{Y})^{i}_{.l}X^{l}-C^{i}_{.kl}X^{k}Y^{l}.
\qquad (15)
\]
If we introduce the matrix $\mathbf{R}^{D}(X,Y):=  (R^{D}(X,Y))^{i}_{.k}  $,
then (14) takes the form
\[
\mathbf{R}^{D}(X,Y)=X(W_{Y})-Y(W_{X})+W_{X}W_{Y}-W_{Y}W_{X}-W_{[X,Y]}.
\qquad (14^\prime )
\]
If $D$ is a linear connection $\nabla $, then, due to (8), we
have:
\[
(R^{\nabla }(X,Y))^{i}_{.j}=R^{i}_{\hbox{.jkl}}X^{k}Y^{l},\qquad (16)
\]
\[
(T^{\nabla }(X,Y))^{i}=T^{i}_{.kl}X^{k}Y^{l},\qquad (17)
\] where $[1-6]$
\[
 R^{i}_{\hbox{.jkl}}:=-E_{l}(\Gamma ^{i}_{.jk})+E_{k}(\Gamma
^{i}_{.jl})-\Gamma ^{m}_{.jk}\Gamma ^{i}_{.ml}+\Gamma ^{m}_{.jl}\Gamma
^{i}_{.mk}-\Gamma ^{i}_{.jm}C^{m}_{.kl}\qquad (18)
\]
\[
 T^{i}_{.kl}:=-(\Gamma
^{i}_{.kl}-\Gamma ^{i}_{.lk})-C^{i}_{.kl},\qquad (19)
\]
are the components of the curvature and torsion tensors of $\nabla $.

In the next sections we shall look for special bases $\{E_{i^\prime }\}$ in
which the components $W$   of an $S$-derivation $D$ vanish along some or
along all vector fields X. Evidently, for this purpose we shall have to solve
$(6^\prime )$ with respect to A under certain conditions.

\medskip
\medskip
  {\bf III. SPECIAL BASES FOR DERIVATIONS IN A NEIGHBORHOOD}

\medskip
In this section we shall solve the  problems  for  existence, uniqueness and
holonomicity of basis or bases $\{E_{i^\prime }\}$ in  which  the components
of a given $(S-)$derivation vanish in some neighborhood U.

{\bf Proposition1:}
The following three statements are equivalent:

  (a) In $U$ the $S$-derivation $D$ is a flat linear connection.

  $(b) D$ is in $U$ curvature free, i.e. $R^{D}=0$, and $D_{X}|_{X=0}=0$.

  (c) For $D$ in $U$ exists basis $\{E_{i^\prime }\}$ such that $W$  =0 for every X.

{\it Proof:} We shall proof this proposition according to the
implications $(a)\Leftrightarrow (b)\Leftrightarrow (c)\Leftrightarrow (a)$.

If (a) is true, then (8) and (16) take place, from where, evidently, follows (b) as the flatness of $\nabla $ means $R^{i}_{\hbox{.jkl}}=0$ in U.

Let (b) be valid. The existence of $\{E_{i^\prime }\}$, in which $W$  $=0$,
is equivalent to the existence of a matrix $A:=  A^{i}_{i^\prime }  $
transforming $\{E_{i}\}$ into $\{E_{i^\prime }\}$ and such that (see
$(6^\prime ))$
\[
0=W=A^{-1}(W_{X}A+X(A)),\qquad (20)
\]
i.e. A must be a solution of $X(A)=-W_{X}A$ for every vector field X. The
integrability conditions for this equation are $0=[X,Y]A= =X(Y(A))-Y(X(A))$
for all commuting vector fields $X$ and $Y$, i.e. $[X,Y]=0$. Using
$X(A)=-W_{X}A$ it can easily be calculated that
\[
[X,Y]A=-(R^{D}(X,Y)+W_{[X,Y]})A.\qquad (21)
\]
So, if (b) is valid, then, due to (1) and (3), we get $W_{X}|_{X=0}=0$,
therefore the integrability conditions for (20) are satisfied and,
consequently, the above pointed transformation exists, i.e. in $\{E_{i^\prime
}\}$ we have $W$  $=0$. Hence, from (b) follows (c).

Let (c) be fulfilled. If we arbitrary fix some basis $\{E_{i}\}$, then the
existence of $\{E_{i^\prime }\}$, in which $W$  $=0$, is equivalent to the
existence of matrix A transforming $\{E_{i}\}$ into $\{E_{i^\prime }\}$ and
which, due to $(6^\prime )$, is such that (20) is valid. From this equation,
we get $W_{X}=-(X(A))A^{-1}$, i.e.
$(W_{X})^{i}_{.j}=-[X^{k}(E_{k}(A^{i}_{i^\prime }))]A^{i^\prime }_{i}$which,
due to (8), means that the $S$-derivation $D$ is a linear connection with
local components $\Gamma ^{i}_{.jk}=-(E_{k}(A^{i}_{i^\prime }))A^{i^\prime
}_{j}$. If in $(14^\prime )$ we substitute $W_{X}= =-(X(A))A^{-1}$and use
$X(A^{-1})=-A^{-1}(X(A))A^{-1}$, we get $R^{D}=R^{\nabla }=0$, i.e. $D$ is a
flat linear connection. So, from (c) follows (a).\blacksquare

The main consequence from proposition 1 is that the question for the existence and the properties of special bases in which the components of an $S$-derivation along every vector field to be zeros in a neighborhood is equivalent to that question for (flat) linear
connections as it is expressed by

{\bf Corollary 2:} In a neighborhood there exists a basis in which the
components of an $S$-derivation along every vector field vanish iff $R^{D}=0$
and $D_{X}|_{X=0}=0$, or iff $D$ is a flat linear connection (whose
components vanish in the same basis).

The following two propositions are almost evident (see $(6^\prime )$ and (15) respectively):

{\bf Proposition 3:} All local bases in which the components of a flat $S$-derivation along every vector field vanish in a neighborhood are obtained from one another by linear transformations with constant coefficients.

{\bf Proposition 4:} The local bases in which the components of a flat $S$-derivation along every vector field vanish in a neighborhood are all holonomic or all anholonomic iff the $S$-derivation has zero or nonzero torsion respectively in that neighborhood.

{\bf Remark:} A stronger result is that the mentioned bases are all holonomic or anholonomic at a given point iff the torsion vanishes or is not a zero respectively at that point. We consider only the above statement as it is the most widely used one of that kind.

Now we shall make some conclusions concerning {\it linear connections.}

{\bf Corollary 5:} In a neighborhood there exists a basis in which the components of a linear connection vanish if and only if this connection is flat in that neighborhood.

{\bf Remark 1:} If the connection is torsion free, this is an old classical result that can be found, e.g., in \S106 from [6], in [3], p.142, or in [5].

{\bf Remark 2:} In [7] an analogous statement is pointed out in the $U_{4}$gravity theory, which states that the $U_{4}$connection components can "always be transformed to zero with respect to a suitable
anholonomic system in $U_{4}"$. This statement suffers from two defects: firstly, generally it is not valid in a neighborhood unless the $U_{4}$connection is not flat, a condition which is not even mentioned in [7], and, secondly, in [7] one finds a "proof" of the cited statement not in a neighborhood, but at a fixed point, which in fact is  not a real proof but only a hint for it as it is a simple counting of the number of conditions which must be satisfied by some independent quantities.

{\it Proof:} If there is a basis $\{E_{i^\prime }\}$ in which the components
of a linear connection $\nabla $ are $\Gamma ^{i^\prime }_{..j^\prime
k^\prime }=0$ , then, in accordance with (18), we have $R^{i^\prime
}_{..j^\prime k^\prime l^\prime }=0$, i.e. $\nabla $ is flat. On the
opposite, let $\nabla $ be flat. Then the $S$-derivation $D$ defined by
$D_{X}:=\nabla _{X}$ has components $W_{X}=\Gamma _{k}X^{k}$ at every point
(see (8)) and, due to (14) and (16), is also flat as $R^{D}=R^{\nabla }=0$
and besides $D_{X}|_{X=0}=\nabla _{X}|_{X=0}=0$. Consequently, by
corollary 2, there exists a basis $\{E_{i^\prime }\}$ such that $W$  =0 for
every X. But $W$  $=\Gamma _{k^\prime }X^{k^\prime }$, so that $\Gamma
_{k^\prime }=0$, i.e. $\Gamma ^{i^\prime }_{..j^\prime k^\prime
}=0.\blacksquare $

In the proof of proposition 1 we constructed a bases in which the components of a flat $S$-derivation vanish by proving that the equation $X(A)+W_{X}A=0, X$ being arbitrary, is integrable with respect to A, and hence $W$  =0 for $E_{i^\prime }=A^{i}_{i^\prime }E_{i}$. But as $D$ is in this case a (flat) linear connection $\nabla  ($see corollary 2), we have $W_{X}=\Gamma _{k}X^{k}$for every X. So, A is a solution of $\Gamma _{k}A+E_{k}(A)=0$, which doesn't depend on $X$, or (see $(6^\prime )) \Gamma $  $=A^{-1}(\Gamma _{k}A+E_{k}(A))=0$, i.e. in $\{E_{i^\prime }\}$ the components of $\nabla $ vanish. In this way one can construct a basis, which generally is anholonomic, and in which the components of a flat linear connection vanish.

{\bf Corollary 6:} All local bases in which the components of a flat linear connection vanish in a neighborhood are obtained from one another by linear transformations with constant coefficients.

  {\it Proof:} The result follows from (8) and proposition 3.\blacksquare

{\bf Corollary 7:} If in a neighborhood there exists a holonomic local basis in which the components of a linear connection vanish, then this connection is torsion free in that neighborhood. On the opposite, if a flat linear connection is torsion free in some neighborhood, then all bases in this neighborhood in which the connection's components vanish are holonomic in it.

{\bf Remark:} In different modifications this result can be found, for instance, in $[2, 3, 5, 6]$.

{\it Proof:} This result is a consequence from $(8), (17), (18)$ and proposition 4.\blacksquare

The question for the existence and the properties of basis (bases) in which
the components of an $S$-derivation $D_{X}$ along a {\it fixed} vector field
$X$ vanish in a neighborhood is not so interesting as the considered until
now problem. That is why we shall only sketch briefly the existence of such
bases in the case $X\mid _{x}\neq 0$ for every $x$ from the neighborhood.

If $\{E_{i^\prime }=A^{i}_{i^\prime }E_{i}\}$ is the looked for basis with
the needed property, $W$  $=0$, then its existence, due to $(6^\prime )$, is
equivalent to the existence of $A:=  A^{i}_{i^\prime }  $ satisfying
$W_{X}A+X(A)=0$ for a given X. As $X$ is fixed, the values of A at different
points are connected through the last equation iff these points lie on one
and the same integral curve (path) of X. So, if $\gamma _{y}:J  \to M$ ($J$
being an ${\Bbb R}$ interval) is the integral curve for $X$ passing through
$y\in M$, i.e. $\gamma _{y}(s_{0})=y$ and
$\dot\gamma_{y}(s)=X_{\gamma_{y}(s)}$, $\dot\gamma_{y}$ being the tangent to
$\gamma_{y}$ vector field, for $s\in J$ and a fixed $s_{0}\in J$, then the
equation $W_{X}A+X(A)=0$ along $\gamma _{y}$ reduces to
$(dA/ds)_{\gamma_{y}(s)}=-W_{X}(\gamma _{y}(s))A(\gamma _{y}(s))$. The
general solution of this equation is
\[ A(\gamma_{y}(s))=Y(s,s_{0};-W_{X}\circ \gamma _{y})B(\gamma _{y}),\qquad
(22)
\]
in which $Y=Y(s,s_{0};Z), Z$ being a matrix function of $s$, is the
unique solution of the initial-value problem [8]
\[ \frac{dY}{ds}=ZY,\quad Y_{s=s_{0}}={\Bbb 1} \qquad (23)
\]
 and the nondegenerate matrix $B$ doesn't depend on $s$ Therefore along the
integral curves of $X$ bases exist in which the components of $D_{X}$vanish.
Hence this is true and at  any point at which $X$ is defined.  Due to
$(6^\prime )$, every two such bases $\{E_{i}\}$ and $\{E_{i^\prime }\}$ are
connected by a linear transformation the  matrix A of which is such that
$X(A)=0$.

\medskip
\medskip
  {\bf IV. SPECIAL BASES FOR DERIVATIONS AT A POINT}

\medskip
The purpose of this section is problems analogous to the ones in the previous section to be investigated but in the case describing the behavior of derivations at a given point.

At first we shall consider $S$-derivations with respect to a {\it fixed vector field}, i.e. we shall deal with a {\it fixed} derivation.

{\bf Proposition 8:} At every point $x_{0}\in M$ and for every fixed vector
field $X$ such that $X|_{x_{0}}\neq 0$ there exists a defined in a
neighborhood of $x_{0}$ local basis $\{E_{i^\prime }\}$ in which the
components of a given $S$-derivation $D_{X}$ along $X$ vanish at $x_{0}$, i.e.
$W$  $(x_{0})=0$ in $\{E_{i^\prime }\}$. If $X|_{x_{0}}=0$, then such a
basis exists if in some basis $\{E_{i}\}$ we have $W_{X}(x_{0})=0$.

{\it Proof:} Let $\{x^{\alpha }\}, \alpha ,\beta =1,\ldots  n$ be local
coordinates in a neighborhood of $x_{0}$and $E_{i}:=B^{\alpha }_{i}\partial /
\partial x^{\alpha }$,
$B:= [B^{\alpha }_{i}]  =:[B^{i}_{\alpha }]^{-1}$
being a nondegenerate matrix function. The existence of $\{E_{i^\prime }\}$
is equivalent to the existence of transformation $\{E_{i}\}  \to
\{E_{i^\prime }=A^{i}_{i^\prime }E_{i}\}$ such that $W$  $(x_{0})=0$ which,
as a consequence of $(6^\prime )$, is equivalent to the existence of $A:=
A^{i}_{i^\prime }  $ with the property
\[
W_{X}(x_{0})A(x_{0})+X(A)_{x_{0}}=0.\qquad (24)
\]
If the basis $\{E_{i^\prime }\}$ exists, from (24) follows that if
$X|_{x_{0}}=0$, then $W_{X}(x_{0})=0$. Vice versa, if $X|_{x_{0}}=0$ and
$W_{X}(x_{0})=0$, then (24) is identically satisfied, i.$e$ in this case any
basis has the needed property.

  So, let us below suppose $X|_{x_{0}}\neq 0$.

In this case with respect to A the equation (24) has infinite number of
solutions, a class of which can be formed by putting
\[
A^{j}_{j^\prime }(x)=a^{j}_{j^\prime k}X^{k}(x_{0})-a^{l}_{j^\prime
k}(W_{X}(x_{0}))^{j}_{.l}B^{k}_{\alpha }(x_{0})(x^{\alpha }-x^{\alpha }_{0})+
\]
\[
\qquad +b^{j}_{j^\prime \alpha \beta }(x)(x^{\alpha }-x^{\alpha
}_{0})(x^{\beta }-x^{\beta }_{0}),\qquad (25)
\]
where $  a^{j}_{j^\prime k}X^{k}(x_{0})  $ is a nondegenerate matrix and the $C^{1}$functions
$b^{j}_{j^\prime \alpha \beta }(x)$ and their derivatives are bounded when $x
\to x_{0}$. In fact, from (25), we find
\[
A^{j}_{j^\prime }(x_{0})=a^{j}_{j^\prime k}X^{k}(x_{0}), \quad
E_{k}(A^{j}_{j^\prime })_{x_{0}}=-a^{l}_{j^\prime
k}(W_{X}(x_{0}))^{j}_{.l}\qquad (26)
\]
which convert (24) into
identity.\blacksquare

{\bf Proposition 9:} If for some $S$-derivation $D_{X}$along a fixed vector
field $X$ there is a local basis in which the components of $D_{X}$ vanish at
a given point, then there exist holonomic bases with this property.

{\it Proof:} If $X|_{x_{0}}=0$, then according to the proof of proposition 8
any basis, if any, including the holonomic ones, have the mentioned property.

If $X|_{x_{0}}\neq 0$, then a class of holonomic bases with the needed
property can be found as follows. For the constructed in the proof of
proposition 8 basis $\{E_{i^\prime }\}$, we get
\[
E_{k^\prime }(A^{j}_{j^\prime })_{x_{0}}=A^{k}_{k^\prime
}(x_{0})E_{k}(A^{j}_{j^\prime })_{x_{0}}=-a^{k}_{k^\prime
m}X^{m}(x_{0})a^{l}_{j^\prime k}(W_{X}(x_{0}))^{j}_{.l},
\qquad (27)
\]
Hence, if we take $a^{k}_{k^\prime m}$ to be of the form $a^{k}_{k^\prime
m}=a_{k^\prime }a^{k}_{m}$, with $a_{k^\prime }\neq 0$ and $\det[a^{k}_{m}]
\neq 0,\infty $, we see that the quantities (27) are symmetric with respect
to $k^\prime $ and $j^\prime $. Consequently, choosing appropriately
$B^{k}_{\alpha }$ (e.g.  $B^{k}_{\alpha }=\delta ^{k}_{\alpha })$, there
exist classes of local coordinates $\{y^{i}\}$ and $\{y^{i^\prime }\}$ in a
neighborhood of $x_{0}$ partially fixed by the conditions:
\[
\frac{\partial y^k}{\partial x^\alpha} \Big|_{x_0}
=B^{k}_{\alpha }(x_{0}), \quad i.e.\
E_{k}|_{x_{0}}= \frac{\partial}{\partial y^k} \Big|_{x_{0}},
\]
\[
\frac{\partial y^j}{\partial y^{j'}} \Big|_{x_0}
=A^{j}_{j^\prime }(x_{0})
=a_{j^\prime}a^{j}_{k}X^{k}(x_{0}), \quad i.e.\
E_{k^\prime }|_{x_{0}}=  \frac{\partial}{\partial y^k}\Big|_{x_{0}},
\]
\[
\frac{\partial^2 y^j}{\partial y^{k'} \partial y^{j'}} \Big|_{x_0}
= E_{k^\prime}(A^{j}_{j^\prime })|_{x_{0}}
= -a_{k^\prime }a_{j^\prime}
[a^{k}_{m}X^{m}(x_{0})a^{l}_{k}(W_{X}(x_{0}))^{j}_{.l}].
\]
Evidently,  the coordinate,  and  hence   holonomic,   basis $\{ \partial
/\partial y^{j^\prime }\}$ defined by some of the coordinates $\{y^{j^\prime
}\}$, satisfying  the above conditions, has the needed property.\blacksquare

Let us now turn our attention to $S$-derivations with respect to {\it
arbitrary vector fields}.

{\bf Proposition 10:} An $S$-derivation $D$ is
at some $x_{0}\in M a$ linear connection iff there is a local basis
$\{E_{i^\prime }\}$ in which the components of $D$ along every vector field
vanish at $x_{0}$.

{\bf Remark:} The $S$-derivation $D$ at $x_{0}$ is a
linear connection if for all $X$ and some, and hence any, basis $\{E_{i}\}$
we have $(cf. (8))$
\[
W_{X}(x_{0})=\Gamma_{k}X^{k}(x_{0}),\qquad (28)
\]
 where $\Gamma _{k}$a re constant matrices. This
means (8) to be valid at $x_{0}$, but it may not be true at $x\neq
x_{0}$.

{\it Proof:} Let $\{x^{i}\}$ be local coordinates in a
neighborhood of $x_{0}$and let $D$ be at $x_{0}$ a linear connection, i.e.
$(28)$ to be valid for some $\Gamma _{k}$. We search for a basis
$\{E_{i^\prime }=A^{i}_{i^\prime }\partial /\partial x^{i}\}$ in which $W$
$=0$. Due to $(6^\prime )$ this is equivalent to $\Gamma
_{k}A(x_{0})+\partial A/\partial x^{k}\mid _{x_{o}}=0$. So, if we define
\[
A(y)=B-\Gamma _{k}B(x^{k}(y)-x^{k}(x_{0}))
+B_{kl}(y)(x^{k}(y)-x^{k}(x_{0}))(x^{l}(y)-x^{l}(x_{0})), \qquad(29)
\]
where $B=$const and $B_{kl}$ and their derivatives are bounded functions when
$y \to x_{0}$, we find
\[
A(x_{0})=B, \quad
\partial A/\partial x^{k}\mid _{x_{o}}=-\Gamma _{k}B. \qquad(30)
\]
Hence $\Gamma _{k}A(x_{0})+\partial A/\partial x^{k}\mid _{x_{o}}\equiv
0$ for all A defined by (30), i.e. the bases $\{E_{i^\prime }=A^{i}_{i^\prime
}\partial /\partial x^{i}\}$ with $  A^{i}_{i^\prime }  =A$ have the needed
property.

On the opposite, let in some $\{E_{i^\prime }\}$ be valid $W$  =0 for every
X. Then, fixing a basis $\{E_{i}=A^{i^\prime }_{i}E_{i^\prime }\}$, from
$(6^\prime )$ we get $W_{X}(x_{0})A(x_{0})+X(A)\mid _{x_{0}}=0$, i.e.
$W_{X}(x_{0})=-X(A)\mid _{x_{0}}A^{-1}(x_{0})$, which means that (28) is
satisfied for $\Gamma _{k}=-E_{k}(A)\mid _{x_{0}}A^{-1}(x_{0})$ and
consequently $D$ is at $x_{0}a$ linear connection.\blacksquare

{\bf Proposition 11:} If there exist bases in which the components of an
$S$-derivation along every vector field vanish at a given point, then they
are obtained from one another by linear transformations whose coefficients
are such that the action of the vectors of these bases on them vanish at the
given point.

{\it Proof:} If $\{E_{i}\}$ and $\{E_{i^\prime }\}$ are such bases at the point $x_{0}$, then $W$  $(x_{0})=W_{X}(x_{0})=0$. Therefore, from $(6^\prime )$, we get $X(A)\mid _{x_{0}}=0$ for every $X$, i.$e E_{i}(A)\mid _{x_{0}}=0.\blacksquare $

{\bf Proposition 12:} If for some $S$-derivation $D$ there is a local holonomic basis in which the components of $D$ along every vector field vanish at a point $x_{0}$, then the torsion of $D$ is zero at $x_{0}$. On the opposite, if $D$ is torsion free at $x_{0}$and bases with the mentioned property exist, then all of them are holonomic at $x_{0}$.

{\it Proof:} If $\{E_{i^\prime }\}$ is a basis with the mentioned property, i.e.
$W$  $(x_{0})=0$ for every $X$, then, using (15), we find $T^{D}(E_{i^\prime
},E_{j^\prime })|_{x_{0}}= =-[E_{i^\prime },E_{j^\prime }]|_{x_{0}}$and
consequently $\{E_{i^\prime }\}$ is holonomic at $x_{0}$, i.e.
$[E_{i^\prime },E_{j^\prime }]|_{x_{0}}=0$,
iff $0=T^{D}(X,Y)|_{x_{0}}=X^{i^\prime
}(x_{0})Y^{j^\prime }(x_{0})(T^{D}(E_{i^\prime },E_{j^\prime })|_{x_{0}})
($see proposition 10 and (28)) for every vector fields $X$ and $Y$, which is
equivalent to $T^{D}|_{x_{0}}=0$.

On the opposite, let $T^{D}|_{x_{0}}=0$. We want to prove that any basis
$\{E_{i^\prime }\}$ in which $W$  $(x_{0})=0$ is holonomic at $x_{0}$. The
holonomicity at $x_{0}$means $0=[E_{i^\prime },E_{j^\prime
}]|_{x_{0}}=-A^{k^\prime }_{k}(E_{j^\prime }(A^{k}_{i^\prime })-E_{i^\prime
}(A^{k}_{j^\prime }))E_{k^\prime }|_{x_{0}}$. But (see proposition 1) the
existence of $\{E_{i^\prime }\}$ is equivalent to $W_{X}(x_{0})=\Gamma
_{k}X^{k}(x_{0})$ for every X. These two facts, combined with (15), show that
$(\Gamma _{k})^{i}_{.j}=(\Gamma _{j})^{i}_{.k}$. Using this and $\Gamma
_{k}A(x_{0})+\partial A/\partial x^{k}|_{x_{0}}=0 ($see the proof of
proposition 10), we find $E_{j^\prime }(A^{k}_{i^\prime
})|_{x_{0}}=-A^{j}_{j^\prime }A^{i}_{i^\prime }(\Gamma
_{j})^{k}_{.i}|_{x_{0}}=E_{i^\prime }(A^{k}_{j^\prime })|_{x_{0}}$and
therefore $[E_{i^\prime },E_{j^\prime }]|_{x_{0}}=0$ (see above), i.e.
$\{E_{i^\prime }\}$ is holonomic at $x_{0}.\blacksquare $

Now we shall apply the above-obtained results to the theory of {\it linear connections}.

{\bf Corollary 13:} For every point $x_{0}$and every linear connection $\nabla $ there exist in a neighborhood of $x_{0}$local bases in which the components of $\nabla $ vanish at $x_{0}$.

{\bf Remark:} For torsion free linear connections this result is well known and is valid in holonomic bases (normal coordinates); see, e.g.: $[2], ch$. III, $\S8; [4]$, p. 120; or $[5], \S36$.

{\it Proof:} From  proposition 10, its proof, (28) and (8) follows the
existence of bases $\{E_{i^\prime }\}$ such that $0=W$  $(x_{0})=\Gamma
_{k^\prime }(x_{0})X^{k^\prime }(x_{0})$ for every $X$, where $\Gamma
_{k^\prime }(x_{0})=  \Gamma ^{i^\prime }_{..j^\prime k^\prime }(x_{0})  $ is
the matrix of the components of $\nabla $ in $\{E_{i^\prime }\}$ at $x_{0}$.
Consequently, as $X$ is arbitrary, $\Gamma _{k^\prime }(x_{0})=0$, i.e.
$\Gamma ^{i^\prime }_{..j^\prime k^\prime }(x_{0})=0.\blacksquare $

  One can easily prove the following three corollaries:

{\bf Corollary 14:} The bases in which the components of a linear connection
$\nabla $ vanish at a point $x_{0}$are obtained from one another by linear
transformations the coefficients of which are such that the action of the
vectors of these bases on them vanishes at $x_{0}$.

{\bf Corollary 15:} In a neighborhood of a given point $x_{0}$there exist
holonomic bases in which the components of a linear connection $\nabla $
vanish at $x_{0}$iff the torsion of $\nabla $ vanishes at $x_{0}$.

{\bf Remark:} This is a classical result that can be found, for instance, in
$[2], ch$. III, \S8 or in [4], p. 120. The same is valid and for corollary 16
below.

{\bf Corollary 16:} For a torsion free linear connection in neighborhood of
any point, local coordinates exist (or, equivalently, holonomic bases) in
which its components vanish at that point.

If $\nabla $ is arbitrary linear connection, then, generally, its torsion is
not zero. But if we define a linear connection $^{s}\nabla $ whose components
are the symmetric part of the ones of $\nabla $, then $^{s}\nabla $ is
torsion free. By corollary 15 for $^{s}\nabla $, local holonomic bases exist
in which its components vanish at any given point. Thus we have proved the
known result (see, e.g., $[2], ch$. III, \S8 and [4], p. 120) that if $\nabla
$ isn't torsion free, then there doesn't exist local holonomic basis in which
the components of $\nabla $ vanish at some point, but there exist local
holonomic bases (coordinates, called normal $[2, 4, 5])$ in which the
symmetric part of the components of $\nabla $ vanish at that point.

\medskip
\medskip
  {\bf V. THE EQUIVALENCE PRINCIPLE}

\medskip
Physically the above results are important in connection with the equivalence principle (see, e.g. $[7, 9]$ and the references in
them).

Usually in a local frame (basis) the gravitational field strength is identified with the components of some linear connection which may be with or without torsion (e.g. the Riemannian one in general relativity [9] or the one in Riemann-Cartan space-times [7]). This linear connection must be compatible with the equivalence principle in a sense that there must exist "local" inertial, called also Lorentz, frames of reference (bases) in which the gravity field strength is "locally" transformed to zero. In our terminology this means the existence of special "local" basis (or bases) in which the connection's components vanish "locally". Above we have put the words "local" and "locally" in inverted commas as they are not well defined here, which is usual for the "physical" literature [9], where they often mean "infinitesimal surrounding of a fixed point of space-time" [7]. The strict meaning of "locally" may be at a point, in a neighborhood, along a path (curve) or on some other submanifold of the space-time. As in the present work we have used the first two of these meanings of "locally" we can make the following conclusions:

(1) All gravity theories based on space-times endowed with a linear connection (e.g. the general relativity [9] and the $U_{4}$theory [7]) are compatible with the equivalence principle at any fixed space-time point, i.e. at any point there exist (local) inertial frames, which generally are anholonomic and may be holonomic ones iff the connection is torsion free (as is e.g. the case of general relativity [9]).

(2) Any gravity theory based on space-time endowed with a flat linear connection is compatible with the equivalence principle in some neighborhood of any space-time point, i.e. for every point there exist its neighborhoods in which there exist (local)
inertial frames (basses) which are holonomic iff the connection is torsion free.

(3) In the above cases the equivalence principle is not at all a principle because it is identically satisfied, namely, it is a corollary from the underlying mathematics for the corresponding gravity theories.

(4) The equivalence principle becomes important if one tries to formulate gravity theories on the base of some (class of) derivations. Generally, it will select those theories which are based on linear connections, i.$e$ only those in which it is identically valid.

\medskip
\medskip
  {\bf VI. REMARKS AND GENERALIZATIONS}

\medskip
As we have seen the linear connections are remarkable among all derivations with their property that in a number of considered here sufficiently general cases they are the only derivations for which special bases in which their components vanish exist.

If one tries to construct a gravity theory based, for example, on linear connections, then he needn't to take into account the equivalence principle for it is identically fulfilled.

This formalism seems to be applicable also to other fields, not only to the gravitational one, namely at least to those of them which are described by linear connections. This suggests the idea for extending the aria of validity of the equivalence principle outside the gravity interaction $(cf. [10])$.

It should be possible to generalize this formalism along paths or on some other submanifolds of the space-time, which will be done elsewhere.

\medskip
\medskip
  {\bf ACKNOWLEDGEMENTS}

\medskip
The author is grateful to Professor Stancho Dimiev (Institute of Mathematics of Bulgarian  Academy of Sciences) and to Dr. Sava Manoff (Institute for Nuclear Research and Nuclear Energy of Bulgarian Academy of Sciences) for the valuable comments and stimulating discussions. He thanks Professor N. A. Chernikov (Joint Institute for Nuclear Research, Dubna, Russia) for the interest to the problems put in this work.

\medskip
\medskip
  REFERENCES

\medskip
1. Dubrovin B., S. P. Novikov, A. Fomenko, Modern Geometry, I. Methods and Applications (Springer Verlag).\par
2. Kobayashi S., K. Nomizu, Foundations of Differential Geometry (Interscience Publishers, New York-London, 1963), vol. {\bf I}.\par
3. Schouten J. A., Ricci-Calculus: An Introduction to Tensor Analysis and its Geometrical Applications (Springer Verlag, Berlin-G\"ottingen-Heidelberg, $1954), 2-nd ed$.\par
4. Lovelock D., H. Rund, Tensors, Differential Forms, and Variational Principals (Wiley-Interscience Publication, John Wiley \& Sons, New York-London-Sydney-Toronto, 1975).\par
5. Norden A. P., Spaces with Affine Connection (Nauka, Moscow, 1976) (In Russian), ch.III, $\S\S35, 36$.\par
6. Rashevskii P. K., Riemannian Geometry and Tensor Analysis (Nauka, Moscow, 1967) (In Russian).\par
7. von der Heyde P., Lett. Nuovo Cimento, vol.{\bf 14}, No.$7, 250,
(1975)$.\par
8. Hartman Ph., Ordinary Differential Equations (John Wiley \& Sons, New York-London-Sydney, $1964), ch$.IV, \S1.\par
9. C. W. Misner, Thorne K. S., Wheeler J. A., Gravitation (W. H. Freeman and Company, San Francisco, 1973); S. Weinberg, Gravitation and Cosmology (J. Wiley \& Sons, New York-London-Sydney-Toronto, 1972).\par
10. E. Kapuszik, Kempchinski Ja., Hozhelya A., Generalized Einstein's
equivalence principle, in: Proceedings of $3-d$ seminar "Gravitational Energy
and Gravitational Waves"{\Large ,} Dubna, $18-20$ may 1990 (JINR,
$P2-91-164,$Dubna, 1991) (In Russian).\par

\end{document}